\documentclass[smallextended,referee,envcountsect]{svjour3} 
\smartqed 
\usepackage{graphicx}
\usepackage{amsmath}
\usepackage{color}
\usepackage{amsfonts}
\usepackage{cite}

\newcommand{\argmax}{\mathop{\mathrm{argmax}}}

\journalname{JOTA}

\begin{document}

\title{Performance Bounds with Curvature for Batched Greedy Optimization}\thanks{This work is supported in part by NSF under award CCF-1422658, and by the CSU Information Science and Technology Center (ISTeC). A preliminary version of a subset of this paper was presented at the \emph{2016 American Control Conference}~\cite{Liu2016}.}

%\subtitle{Performance Bounds for Batched Greedy Optimization}

\author{Yajing Liu   \and  Zhenliang~Zhang \and  Edwin~ K.~P.~Chong  \and Ali~Pezeshki }

\institute{Yajing Liu (Corresponding author), Edwin K. P. Chong, Ali Pezeshki \at
           Colorado State University \\
              Fort Collins, CO, USA\\
             yajing.liu@ymail.com, edwin.chong@colostate.edu, ali.pezeshki@colostate.edu
           \and
          Zhenliang Zhang \at
         Intel Labs\\
          Hillsboro, Oregon, USA\\
        zzl.csu@gmail.com
}

\date{Communicated by Weibo Gong\\
Received: date / Accepted: date}
%The correct dates will be entered by the editor.

\maketitle

\begin{abstract}
The batched greedy strategy is an approximation algorithm to maximize a set function subject to a matroid constraint. Starting with the empty set, the batched greedy strategy iteratively adds  to the current solution set a batch of elements that results in the largest gain in the objective function  while satisfying the matroid constraints. In this paper, we develop bounds on the performance of the  batched greedy strategy relative to the optimal strategy in terms of a parameter  called the total  batched curvature. We show that when the objective function is a polymatroid set function, the batched greedy strategy satisfies  a harmonic bound   for a general matroid constraint and   an exponential bound  for a uniform matroid constraint, both in terms of the total batched curvature.
We also study the behavior of the bounds as functions of the batch size. Specifically, we prove that the harmonic bound for a general matroid is nondecreasing in the batch size and the exponential bound for a uniform matroid is nondecreasing in the batch size under the condition that the batch size divides the rank of the uniform matroid. Finally, we illustrate our results by considering a task scheduling problem and an adaptive sensing problem.
\end{abstract}
\keywords{Curvature\and Greedy \and Matroid \and Polymatroid \and Submodular}
\subclass{90C27 \and 90C59}

%All acknowledgements should be placed in the back of the paper after Conclusions..

\section{Introduction}

A variety of combinatorial optimization problems such as generalized assignment (see, e.g., \cite{streeter2008online, FeigeVondrak, Fleischer2006, Cohen2006,  Shmoys1993}), welfare maximization (see, e.g., \cite{Vondrak2011, Korula2015, Vondrak2008}), maximum coverage (see, e.g., \cite{K-cover1998, Feige1998, Khuller1999}), maximal covering location  (see, e.g., \cite{Fisher1977,  Location, Church1974,  Pirkul1991}), and sensor placement (see, e.g., \cite{LiC12, SensorPlacement, Chen2005, Ragi2015}) can be formulated as a problem of maximizing a set function subject to a matroid constraint. More precisely,  the objective function maps the power set of a ground set to real numbers, and the constraint is that any feasible set is from a non-empty collection of subsets of the ground set satisfying matroid constraints.

Finding the optimal solution to the problem above in general is NP-hard. The \emph{greedy strategy} provides a computationally feasible approach, which starts with the empty set, and then iteratively adds to the current solution set one element  that results in  the largest gain in the objective function, while satisfying the matroid constraints. This scheme is a special case of the \emph{batched greedy strategy} with batch size equal to 1. For general batch size (greater than 1), the batched greedy strategy  starts with the empty set but iteratively adds, to the current
solution set, a batch of elements with the largest gain in the objective function under the constraints. A detailed definition of the batched greedy strategy is given in Section~\ref{kbatchgreedy}. 

The performance of the batched greedy strategy {with batch size equal to 1}  has been extensively investigated in \cite{nemhauser19781,  nemhauser1978, hausmann1980, conforti1984submodular, vondrak2010submodularity, sviridenko2015}, which we will review in Section~\ref{previouswork}.
The performance of the batched greedy strategy for general {batch size}, however,  has received little attention, notable exceptions being Nemhauser \emph{et~al.}~\cite{nemhauser19781} and Hausmann \emph{et al.}~\cite{hausmann1980}, which we will review in Section~\ref{previouswork}.  Although Nemhauser \emph{et al.}~\cite{nemhauser19781} and Hausmann \emph{et al.}~\cite{hausmann1980} investigated the performance of the batched  greedy strategy, they only considered uniform matroid constraints and independence system constraints, respectively. This prompts us to investigate the performance of the batched strategy more comprehensively. In particular, it is of interest to extend the notion of total curvature to the batched case.

Our main contribution is that we derive bounds in terms of the total batched curvature for the performance of the batched greedy strategy under  general matroid and  uniform matroid constraints. We also study the behavior of the bounds as functions of the batch size, by comparing the values of the total batched curvature for different batch sizes and investigating the monotoneity of the bounds. It is \emph{not} our claim that we are proposing a new algorithm (the batched greedy strategy) or even that we are advocating the use of such an algorithm. Our contribution is to provide bounds on the performance of the batched greedy strategy, which we consider to be a rather natural extension of the greedy strategy. As we argue below and as the reader can see in the remainder of this paper, going from the case of batch size equal to 1 to the general case (batch size greater than 1) is highly nontrivial.

In \cite{conforti1984submodular}, Conforti and Cornu{\'e}jols  provided  performance bounds for the greedy strategy in terms of the total curvature under general matroid constraints and uniform matroid constraints. It might be tempting to think that bounds for the batched case can be derived in a straightforward way from the results of batch size equal to 1 by \emph{lifting}, which is to treat each batch-sized set of elements
chosen by the batched greedy strategy as a single action,
and then appeal to the results for the case of batch size equal to 1. However, it turns out that lifting \emph{does not work} for a general batched greedy strategy (batch size greater than 1) for the following two reasons. First, the collection of sets created by satisfying the batched greedy strategy is not a matroid in general; we will demonstrate this by an example in Appendix~\ref{AppendixA}. Second, the last step of the batched greedy strategy may select elements with a number less than the batch size, because the cardinality of the maximal set in the matroid may not be divisible by the batch size. For the cases considered in \cite{Liu2016}, lifting fails to work for the first reason.

The batched greedy strategy requires an exponential number of evaluations of the objective function if using exhaustive search. When the batch size is equal to the cardinality of the maximal set in the matroid, the batched greedy strategy coincides with the optimal strategy. It might be tempting to expect that the batched strategy with batch size  greater than 1 outperforms the usual greedy strategy, albeit at the expense of increasing computational complexity. Indeed, the Monte Carlo simulations performed in  \cite{Chandu2015}  for the  \emph{maximum coverage} problem  show that the batched greedy strategy {with batch size greater than 1} provides better approximation than the usual greedy strategy in many cases. However, it is also evident from their simulation that this is not always the case.  In Appendix~\ref{AppendixB}, we provide two examples of the {maximum coverage problem} where the usual greedy strategy performs better than the batched greedy strategy with batch size 2.

In Section~\ref{sc:II}, we first introduce some definitions and review the previous results. Then, we review Lemmas~1.1 and~1.2 from \cite{vondrak2010}, which we will use to derive performance bounds for the batched greedy strategy under a uniform matroid constraint. In Section~\ref{sc:III}, we define the total batched curvature and then we provide a harmonic bound and an exponential bound for the batched greedy strategy under a general matroid constraint and a uniform matroid constraint, respectively, both in terms of the total batched curvature. When the total batched curvature equals 1, the bound for a uniform matroid becomes the bound in \cite{nemhauser19781}. When  the batch size equals 1, the bounds reduce to the bounds derived in  \cite{conforti1984submodular}. When the batch size divides the rank of the uniform matroid, the bounds become the bounds in \cite{Liu2016}. We also prove that the batched curvature is nonincreasing in the batch size when the objective function is  a polymatroid set function. This implies that the larger the batch size, the better the harmonic bound for a general matroid and when  the batch size divides  the rank of the uniform matroid, the larger the batch size, the better the exponential bound for a uniform matroid. In Section~\ref{sc:IV}, we present a task scheduling problem and an adaptive sensing problem to demonstrate our results. 

\section{Preliminaries}\label{sc:II}
\subsection{Polymatroid Set Functions and Curvature}
The definitions and terminology in this paragraph are standard (see, e.g., \cite{Edmonds, Boros2003, Tutte}), but are included for completeness.
Let $X$ be a finite set, and $\mathcal{I}$ be a non-empty collection of subsets of  $X$. Given a pair $(X,\mathcal{I})$, the collection $\mathcal{I}$ is said to be \emph{hereditary} iff it satisfies property i below and  has the \emph{augmentation} property iff it satisfies property ii below:
\begin{itemize}
\item  [i.] For all $B\in\mathcal{I}$, any set $A\subseteq B$ is also in $\mathcal{I}$.
\item  [ii.] For any $A,B\in \mathcal{I}$, if $|B|>|A|$, then there exists $j\in B\setminus A$ such that $A\cup\{j\}\in\mathcal{I}$.
\end{itemize}
The pair $(X,\mathcal{I})$ is called a \emph{matroid} iff it satisfies both properties i and ii. The pair $(X,\mathcal{I})$ is called a \emph{uniform matroid} iff $\mathcal{I}=\{S\subseteq X: |S|\leq K\}$ for a given $K$, called the \emph{rank} of  $(X,\mathcal{I})$.

Let $2^X$ denote the power set of $X$, and define a set function $f$: $2^X\longrightarrow {\rm I\!R}$.
The set function $f$ is said to be \emph{nondecreasing} and \emph{submodular} iff it satisfies properties~1 and~2 below, respectively:
\begin{itemize}
\item [1.] For any $A\subseteq B\subseteq X$, $f(A)\leq f(B)$.
\item [2.] For any $A\subseteq B\subseteq X$ and $j\in X\setminus B$, $f(A\cup\{j\})-f(A)\geq f(B\cup\{j\})-f(B)$.
\end{itemize}
A set function $f$: $2^X\longrightarrow {\rm I\!R}$ is called a \emph{polymatroid set function} iff it is nondecreasing, submodular, and $f(\emptyset)=0$, where $\emptyset$ denotes the empty set.
The submodularity in property~2 means that the additional value accruing from an extra action decreases as the size of the input set increases, and is also called the \emph{diminishing-return} property in economics.
Submodularity implies that for any $A\subseteq B\subseteq X$ and $T\subseteq X\setminus B$, 
\begin{equation}
\label{eqn:submodularimplies}
f(A\cup T)-f(A)\geq f(B\cup T)-f(B).
\end{equation}
For convenience, we denote the incremental value of adding a set $T$ to the set $A\subseteq X$ as $\varrho_T(A)=f(A\cup T)-f(A)$ (following the notation of \cite{conforti1984submodular}).

The \emph{total curvature} of a set function $f$ is defined as \cite{conforti1984submodular}
$$c:=\max_{j\in X^*}\left\{1-\frac{\varrho_j({X\setminus\{j\}})}{\varrho_j(\emptyset)}\right\},$$
where $X^*=\{j\in X: \varrho_j(\emptyset)\neq 0\}$. Note that $0\leq c\leq 1$ when $f$ is a polymatroid set function, and $c=0$ if and only if $f$ is \emph{additive}, i.e., for any set $A\subseteq X$, $f(A)=\sum_{i\in A}f(\{i\})$. When $c=0$, it is easy to check that the greedy strategy coincides with the optimal strategy. So in the rest of the paper, when we assume that $f$ is a polymatroid set function, we only consider $c\in]0,1]$.
\subsection{Review of Previous Work}
\label{previouswork}
Before we review the previous work, we formulate the optimization problem formally as follows:
\begin{align}\label{eqn:1}
\begin{array}{l}
\text{maximize} \ \    f(M), \ \quad \text{subject to} \ \  M\in \mathcal{I},
\end{array}
\end{align}
where  $\mathcal{I}$ is a non-empty collection of subsets of a finite set $X$, and $f$ is a real-valued set function defined on the power set $2^X$ of $X$.

For convenience, in the rest of the paper we will use \emph{$k$-batch greedy strategy} to denote the batched greedy strategy with batch size $k$. So, the $1$-batch greedy strategy denotes the usual greedy strategy.

Nemhauser \emph{et al.}~\cite{nemhauser19781, nemhauser1978} proved that, when $f$ is a polymatroid set function, the $1$-batch greedy strategy yields a $1/2$-approximation\footnotemark
% ...
\footnotetext{The term $\beta$-approximation means that $f(G)/f(O)\geq \beta$, where $G$ and $O$ denote a greedy solution and an optimal solution, respectively.} for a general matroid and a $(1-e^{-1})$-approximation for a uniform matroid. By introducing the total curvature $c$,
Conforti and Cornu{\'e}jols~\cite{conforti1984submodular} showed that, when $f$ is  a polymatroid set function, the $1$-batch greedy strategy achieves   a $1/(1+c)$-approximation for a general matroid and a $(1-e^{-c})/{c}$-approximation for a uniform matroid. For a polymatroid  set function $f$,  the total curvature $c$ takes values on the interval $ ]0,1]$. In this case, we have $1/(1+c)\geq1/2$ and $(1-e^{-c})/c\geq (1-e^{-1})$, which implies that the bounds $1/(1+c)$ and $(1-e^{-c})/c$ are stronger than
the bounds $1/2$ and $(1-e^{-1})$ in \cite{nemhauser1978} and \cite{nemhauser19781}, respectively. Vondr{\'a}k~\cite{vondrak2010submodularity} proved that, when $f$ is  a polymatroid set function, the continuous greedy strategy gives a $(1-e^{-c})/c$-approximation for any matroid. Sviridenko \emph{et al.}~\cite{sviridenko2015} proved that, when  $f$ is  a polymatroid set function, a modified continuous greedy strategy gives a $(1-ce^{-1})$-approximation for any matroid, the first improvement over the greedy $(1-e^{-c})/{c}$-approximation of Conforti and Cornu{\'e}jols from~\cite{conforti1984submodular}.

Nemhauser \emph{et al.}~\cite{nemhauser19781} proved that, when $f$ is  a polymatroid set function and $(X,\mathcal{I})$ is a uniform matroid of rank $K=kl+m$ ($l$ and $m$ are nonnegative integers and $0< m\leq k$), the $k$-batch greedy strategy achieves a $\gamma$-approximation, where $\gamma=(1-(1-\frac{m}{k(l+1)})(1-\frac{1}{l+1})^l)$.  Hausmann \emph{et al.}~\cite{hausmann1980} showed that, when $f$ is  a polymatroid set function and $(X,\mathcal{I})$ is an independence system, the $k$-batch greedy strategy achieves  a $q(X,\mathcal{I})$-approximation, where $q(X,\mathcal{I})$ is the rank quotient defined in \cite{hausmann1980}.

\subsection{Performance Bounds in Terms of  Total Curvature}
\label{reviewbounds}
In this section, we review two theorems from \cite{conforti1984submodular}, which bound the performance of the $1$-batch greedy strategy using the total curvature $c$ for  general matroid constraints and  uniform matroid constraints. These bounds are special cases of the bounds we derive in Section~3.2 for $k=1$. 

We first define optimal and greedy solutions for problem~\eqref{eqn:1} as follows:

\emph{Optimal solution}: 
Consider problem \eqref{eqn:1} of finding a set that maximizes $f$ under the constraint $M\in\mathcal{I}$. 
We call a solution of this problem an \emph{optimal solution} and
denote it by $O$, i.e., 
\begin{align*}
O\in\mathop{\argmax}\limits_{M\in \mathcal{I}} f(M),
\end{align*}
 where $\argmax$ denotes the set of actions that maximize $ f(\cdot)$. 

\emph{$1$-batch greedy solution}:
A set $G=\{g_1,g_2,\ldots,g_{k}\}$ is called a \emph{$1$-batch greedy} solution  if
\begin{align*}
g_1&\in\mathop{\argmax}\limits_{g\in{X}} f(\{g\}),
\end{align*} 
and 
for $i=2,\ldots,k,$
\begin{align*}
g_i&\in\mathop{\argmax}\limits_{g\in{X}} f(\{g_1,g_2,\ldots,g_{i-1},g\}).
\end{align*}

\begin{theorem}\textup{{\cite{conforti1984submodular}}}
\label{Theorem2.1}
Let  $(X,\mathcal{I})$ be a matroid and $f$: $2^X\longrightarrow~{\rm I\!R}$ be a polymatroid set function with total curvature $c$. Then, any $1$-batch greedy solution $G$ satisfies
$$\frac{f(G)}{f(O)}\geq \frac{1}{1+c},$$
where $O$ is any optimal solution to problem \textup{(\ref{eqn:1})}.
\end{theorem}

 When $f$ is  a polymatroid set function, we have $c\in ]0,1]$, and therefore $1/(1+c)\in[1/2,1[$. Theorem \ref{Theorem2.1} applies to any matroid. This means that the bound ${1}/(1+c)$ holds for a uniform matroid  too. Theorem \ref{Theorem2.2} below provides a tighter bound when $(X,\mathcal{I})$ is a uniform matroid.

\begin{theorem}\textup{\cite{conforti1984submodular}}
\label{Theorem2.2}
\sloppy{Let $(X,\mathcal{I})$ be a uniform matroid of rank $K$. Further, let $f:2^X\longrightarrow{\rm I\!R}$ be a polymatroid set function with total curvature $c$. Then, any $1$-batch greedy solution $G$ satisfies
\begin{align*}
\frac{f(G)}{f(O)}&\geq\frac{1}{c}\left(1-\left(1-\frac{c}{K}\right)^K\right)\geq \frac{1}{c}\left(1-e^{-c}\right),
\end{align*}
where $O$ is any optimal solution to problem \textup{(\ref{eqn:1})}.}
\end{theorem}
\vspace{2mm}

The function $(1-e^{-c})/c$ is a nonincreasing function of $c$, and therefore $(1-e^{-c})/c\in[1-e^{-1},1[$ when $f$ is a polymatroid set function. Also it is easy to check that $(1-e^{-c})/{c}\geq 1/(1+c)$ for $c\in]0,1]$, which implies that the bound $(1-e^{-c})/{c}$ is stronger than the bound $1/(1+c)$ in Theorem \ref{Theorem2.1}.
\subsection{Properties of Submodular Functions}
\label{reviewtheorems}
The following two lemmas from \cite{vondrak2010}, stating some technical properties of submodular functions, will be useful to derive performance bounds for the $k$-batch greedy strategy under a uniform matroid constraint.
\begin{lemma}\textup{\cite{vondrak2010}}
\label{lemma1}
Let  $f$: $2^X\longrightarrow {\rm I\!R}$ be a submodular set function. Given $A, B\subseteq X$, let $\{M_1,\ldots, M_r\}$ be a collection of subsets of $B\setminus A$ such that each element of $B\setminus A$ appears in exactly $p$ of these subsets. Then,
\[\sum\limits_{i=1}^r\varrho_{M_i}(A)\geq p\varrho_B(A).\]
\end{lemma}
\vspace{1mm}
\begin{lemma}\textup{\cite{vondrak2010}}
\label{lemma2}
Let  $f$: $2^X\longrightarrow {\rm I\!R}$ be a submodular set function. Given $A'\subseteq A\subseteq X$, let $\{T_1,\ldots, T_s\}$ be a collection of subsets of $A\setminus A'$ such that each element of $A\setminus A'$ appears in exactly $q$ of these subsets. Then, 
\[\sum\limits_{i=1}^s\varrho_{T_i}(A\setminus T_i)\leq q\varrho_{A\setminus A'}(A').\]
\end{lemma}

\section{Main Results}\label{sc:III}

In this section, first we  define the $k$-batch greedy strategy and the total $k$-batch curvature $c_k$ that will be used for deriving harmonic and exponential bounds. Then we derive performance bounds for the $k$-batch greedy strategy in terms of  $c_k$ under general matroid constraints and under uniform matroid constraints. Moreover, we  study the behavior of the bounds as functions of the batch size $k$.
\subsection{$k$-Batch Greedy Strategy}
\label{kbatchgreedy}
We write the cardinality of the  maximal set in $\mathcal{I}$ as $K=kl+m$, where $l,m$ are nonnegative integers and $0<m\leq k$. Note that $m$ is not necessarily the remainder of $K/k$, because $m$ could be equal to $k$. This happens when $k$ divides $K$. The $k$-batch greedy strategy is as follows:

Step~1: Let $S^0=\emptyset$ and $t=0$.

Step~2: Select $J_{t+1}\subseteq X\setminus S^t$ for which $|J_{t+1}|=k$, $S^t\cup J_{t+1}\in\mathcal{I}$, and 
\begin{align*}
f(S^t\cup J_{t+1})=\max\limits_{J\subseteq X\setminus S^t\ \text{and}\ |J|=k }f(S^t\cup J);
\end{align*}
 then set $S^{t+1}=S^t\cup J_{t+1}$.

Step~3:  If $t+1< l$, set $t=t+1$, and repeat Step~2.

Step~4: If $t+1=l$, select $J_{l+1}\subseteq X\setminus S^l$ such that $|J_{l+1}|=m$, $S^l\cup J_{l+1}\in\mathcal{I}$, and 
\[f(S^l\cup J_{l+1})=\max\limits_{J\subseteq X\setminus S^l\ \text{and}\ |J|=m}f(S^l\cup J).\]

Step~5: Return the set $S=S^l\cup J_{l+1}$ and terminate.

Any set generated by the above procedure is called a \emph{$k$-batch greedy solution}.

The difference between a $k$-batch greedy strategy for a general matroid and that for a uniform matroid is that at each step $t$ ($0\leq t\leq l$), we have to check whether $J_{t+1}\subseteq X\setminus S^t$ satisfies $S^t\cup J_{t+1}\in\mathcal{I}$ for a general matroid while $S^t\cup J_{t+1}\in\mathcal{I}$  always holds for a uniform matroid.

\subsection{Performance Bounds in Terms of  Total $k$-Batch Curvature}
Similar to the definition of the total curvature $c$ in \cite{conforti1984submodular}, we define the \emph{total $k$-batch curvature} $c_k$ for a given $k$ as
\begin{equation}
\label{totalkbatchcurvature}
c_k:=\max\limits_{I\in \hat{X}}\left\{1-\frac{\varrho_I(X\setminus I)}{\varrho_I(\emptyset)}\right\},
\end{equation}
 where $\hat{X}=\{I\subseteq X: \varrho_I(\emptyset)\neq 0 \ \text{and}\ |I|=k\}$.

 The following  proposition will be applied to derive our bounds in terms of $c_k$ for both general matroid constraints and  uniform matroid constraints.
\begin{proposition}
\label{Prop1}
If $f:2^X\longrightarrow {\rm I\!R}$ is a submodular set function, $A,B\subseteq X$, and $\{M_1,\ldots, M_r\}$ is a partition of $B\setminus A$, then 
\begin{equation}
\label{ineq:prop1}
f(A\cup B)\leq f(A)+\sum\limits_{i:M_i\subseteq B\setminus A}\varrho_{M_i}(A).
\end{equation}
\end{proposition}
\vspace{1mm}
{\it Proof}
By the assumption that $\{M_1,\ldots, M_r\}$ is a partition of $B\setminus A$ and by submodularity (see  inequality (\ref{eqn:submodularimplies})),  we have 
 \begin{align*}
 f(A\cup B)-f(A)&=f(A\cup \bigcup_{j=1}^r M_j)-f(A)\\
 &=\sum\limits_{i=1}^r \varrho_{M_i}(A\cup\bigcup_{j=1}^{i-1}M_j)\leq \sum\limits_{i:M_i\subseteq B\setminus A}\varrho_{M_i}(A),
 \end{align*}
which implies  inequality (\ref{ineq:prop1}).
\qed
The following proposition in terms of the total $k$-batch curvature $c_k$ will be applied to derive our  bounds under  general matroid constraints.
 \begin{proposition}
 \label{Prop3}
Let $f:2^X\longrightarrow {\rm I\!R}$ be a polymatroid set function. Given a set $B\subseteq X$,  a sequence of $t$ $ (t>0 )$ sets $A^i=\bigcup_{j=1}^iI_j$ with $I_j\subseteq X$ and $|I_j|=k$ for $1\leq j\leq t$, and a partition $\{M_1,\ldots, M_r\}$  of $B\setminus A^t$, we have 
 \begin{equation}
 \label{ineq:prop3}
  f(B)\leq {c}_k\sum\limits_{i:I_i\subseteq A^t\setminus B}\varrho_{I_i}(A^{i-1})+\sum\limits_{i:I_i\subseteq B\cap A^t}\varrho_{I_i}(A^{i-1})+\sum\limits_{i:M_i\subseteq B\setminus A^t}\varrho_{M_i}(A^t).
 \end{equation}
 \end{proposition}
{\it Proof}
By the definition of $A^t$, we write 
 \begin{align*}
 f(A^t\cup B)-f(B)=\sum\limits_{i=1}^t\varrho_{I_i}(B\cup A^{i-1})=\sum\limits_{i:I_i\subseteq A^t\setminus B}\varrho_{I_i}(B\cup A^{i-1}).
 \end{align*}
 By submodularity (see  inequality (\ref{eqn:submodularimplies})), we have 
\begin{equation}
\label{rhoineq1}
\varrho_{I_i}(B\cup A^{i-1})\geq \varrho_{I_i}(X\setminus I_i)
\end{equation}
and 
\begin{equation}
\label{rhoineq2}
\varrho_{I_i}(\emptyset)\geq \varrho_{I_i}(A^{i-1})
\end{equation}
for $1\leq i\leq t$.
By the definition of the total $k$-batch curvature $c_k$, we have
\[1-\frac{\varrho_{I_i}(X\setminus I_i)}{\varrho_{I_i}(\emptyset)}\leq c_k\]
 for $1\leq i\leq t$,
which implies that 
\[\varrho_{I_i}(X\setminus I_i)\geq (1-c_k)\varrho_{I_i}(\emptyset).\]
Combining the above inequality with (\ref{rhoineq1}) and (\ref{rhoineq2}), we have 
\[\varrho_{I_i}(B\cup A^{i-1})\geq \varrho_{I_i}(X\setminus I_i)\geq (1-c_k)\varrho_{I_i}(\emptyset)\geq (1-c_k)\varrho_{I_i}(A^{i-1})\]
for $1\leq i\leq t$.
Using the above inequality, we have 
\begin{align}
\label{ineqprop61}
 f(A^t\cup B)-f(B)&=\sum\limits_{i:I_i\subseteq A^t\setminus B}\varrho_{I_i}(B\cup A^{i-1})\nonumber\\
 &\geq (1-{c}_k)\sum\limits_{i:I_i\subseteq A^t\setminus B}\varrho_{I_i}( A^{i-1}).
 \end{align}
 By Proposition~\ref{Prop1}, we have 
\begin{equation}
\label{ineqprop62}
 f(A^t\cup B)\leq f(A^t)+\sum\limits_{i:M_i\subseteq B\setminus A^t}\varrho_{M_i}(A^t).
\end{equation}
Combining inequalities~(\ref{ineqprop61}) and~(\ref{ineqprop62}) results in 
\[f(B)\leq f(A^t)+\sum\limits_{i:M_i\subseteq B\setminus A^t}\varrho_{M_i}(A^t)-(1-{c}_k)\sum\limits_{i:I_i\subseteq A^t\setminus B}\varrho_{I_i}( A^{i-1}).\]
Substituting $f(A^t)$ into the above inequality by the identity \[f(A^t)=\sum\limits_{i:I_i\subseteq A^t\setminus B}\varrho_{I_i}( A^{i-1})+\sum\limits_{i:I_i\subseteq B\cap A^t}\varrho_{I_i}(A^{i-1}),\] we  get  inequality (\ref{ineq:prop3}).
 \qed

Recall that in Section~3.1, we defined  $J_i$ as the set selected by the $k$-batch greedy strategy at stage $i$  and $S^i=\bigcup_{j=1}^{i}J_j$ as the set selected by the $k$-batch greedy strategy for the first $i$ stages, where  $1\leq i\leq l+1$, $|J_i|=k$ for $1\leq i\leq l$,  $|J_{l+1}|=m$, and $K=kl+m$ with $l\geq 0$ and $0< m\leq k$ being integers. When the pair $(X,\mathcal{I})$ is a  matroid,  by  the augmentation property of a  matroid and  the previous assumption that the maximal cardinality of $\mathcal{I}$ is $K$, we have that any optimal solution can be augmented to a set of length $K$. Assume that  $O=\{o_1,\ldots, o_K\}$ is an optimal solution to problem~\eqref{eqn:1}. Let $S=S^{l+1}$ be a $k$-batch greedy solution. We now state and prove the following lemma, which will be used to derive the harmonic bound for general matroid constraints in Theorem~\ref{Theorem3.3}.

\begin{lemma}
\label{lemma3}
Let $S$ be a $k$-batch greedy solution and $O=\{o_1,\ldots, o_K\}$ be an optimal solution. Then the following statements hold:
\begin{itemize}
\item [a.] There exists a partition $\{J_i'\}_{i=1}^{l+1}$ of $O$ with $|J_i'|=k$ for $1\leq i\leq l$ and $|J_{l+1}'|=m$  such that $\varrho_{J_i'}(S^{i-1})\leq \varrho_{J_i}(S^{i-1})$. Furthermore, if $J_i'\subseteq O\cap S$, then $J_i'=J_i$.

\item [b.]  If $J_i'\subseteq O\setminus S^l$ for $1\leq i\leq l$, then $J_i\subseteq S^l\setminus O$.
\end{itemize}
\end{lemma} 
{\it Proof}
We begin by proving a. First, we prove that there exists $J_{l+1}'\subseteq O\setminus S^{l}$ such that $S^l\cup J_{l+1}'\in\mathcal{I}$ and  $\varrho_{J_{l+1}'}(S^l)\leq \varrho_{J_{l+1}}(S^l)$. 
 By definition, $|O|=K$ and $S^l=kl=K-m$. Using the augmentation property, there exists one element $o_{i_1}\in O\setminus S^l$ such that $S^l\cup\{o_{i_1}\}\in\mathcal{I}$. Consider $S^l\cup\{o_{i_1}\}$ and $O$. Using the augmentation property again, there exists one element $o_{i_2}\in O\setminus S^l\setminus \{o_{i_1}\}$ such that $S^l\cup \{o_{i_1}, o_{i_2}\}\in \mathcal{I}$. Using the augmentation property $(m-2)$ more times, we have that there exists $J_{l+1}'=\{o_{i_1},\ldots, o_{i_m}\}\subseteq O\setminus S^l$ such that $S^l\cup J_{l+1}'\in\mathcal{I}$. By the $k$-batch greedy strategy, we have $\varrho_{J_{l+1}'}(S^l)\leq \varrho_{J_{l+1}}(S^l)$.  If $J_{l+1}\subseteq O$, we can set $J_{l+1}'=J_{l+1}$.

Then similar to the proof in \cite{nemhauser19781}, we will prove statement a by backward induction on $i$ for $i=l,l-1,\ldots, 1$. 
Assume that $J_i'$ satisfies the inequality $\varrho_{J_i'}(S^{i-1})\leq \varrho_{J_i}(S^{i-1})$ for $i>j$, and let $O^j=O\setminus \bigcup_{i>j} J_i'$. Consider the sets $S^{j-1}$ and $O^j$. By definition, $|S^{j-1}|=(j-1) k$ and $|O^j|=j k$.  Using the augmentation property, we have that there exists one element $o_{j_1}~\in~O^j~\setminus~S^{j-1}$ such that $S^{j-1}\cup\{o_{j_1}\}\in\mathcal{I}$. Next consider $S^{j-1}\cup\{o_{j_1}\}$ and $O^j$. Using the augmentation property again, there exists one element $o_{j_2}\in O^j\setminus S^{j-1}\setminus\{o_{j_1}\}$ such that $S^{j-1}\cup\{o_{j_1}, o_{j_2}\}\in\mathcal{I}$. Similar to the process above, using the augmentation property $(k-2)$ more times, finally we have that there exists $J_j'=\{o_{j_1},\ldots,o_{j_k}\}\subseteq O^j\setminus S^{j-1}$ such that $S^{j-1}\cup J_j'\in \mathcal{I}$. By the $k$-batch greedy strategy, we have that $\varrho_{J_j'}(S^{j-1})\leq \varrho_{J_j}(S^{j-1})$. Furthermore, if $J_j\subseteq O^j$, we can set $J_j'=J_j$. This completes the proof of statement a.

Now we prove statement b by contradiction. Consider the negation of statement b, i.e.,  if $J_i'\subseteq O\setminus S^l$ for $1\leq i\leq l$, then $J_i\subseteq O$. By the argument in the second paragraph of the proof of statement~a, we have that if $J_i\subseteq O$ for $1\leq i\leq l$, then $J_i=J_i'$. By the assumption that $J_i'\subseteq O\setminus S^l$ for $1\leq i\leq l$, we have $J_i\subseteq O\setminus S^l$ for $1\leq i\leq l$, which contradicts  the fact that $J_i\subseteq S^l$ for $1\leq i\leq l$. 
This completes the proof of statement b.
\qed

The following theorem  presents our performance bound in terms of the total $k$-batch curvature $c_k$ for the $k$-batch greedy strategy under a general matroid.
\begin{theorem}
\label{Theorem3.3}
Let  $(X,\mathcal{I})$ be a general matroid and $f:2^X\longrightarrow {\rm I\!R}$ be a polymatroid set function. Then, any $k$-batch greedy solution $S$ satisfies 
\begin{equation}
\label{ineq:generalbound}
\frac{f(S)}{f(O)}\geq \frac{1}{1+c_k}.
\end{equation}
\end{theorem}
\vspace{2mm}
{\it Proof}
Let $\{P_{i_1},\ldots, P_{i_r}\}$ be a partition of $O\setminus S^l$ satisfying that $P_{i_j}\subseteq J_{i_j}'$ for $1\leq j\leq r$. The way to find $\{P_{i_1},\ldots, P_{i_r}\}$ is as follows:
first list all of the actions in $O\setminus S^l$, then let $P_i$ be its subset consisting of  actions belonging to $J_i'$, i.e., $P_i=(O\setminus S^l)\cap J_i'$. Finally, extract the nonempty sets from $\{P_i\}_{i=1}^{l+1}$ as $\{P_{i_1},\ldots, P_{i_r}\}.$

Recall that $S^i=\cup_{j=1}^i J_j$ for $1\leq i\leq l$ as defined in Section~3.1.  Then using  Proposition~\ref{Prop3}, with $A^t=S^l$ and $S=O$  results in
 \begin{align}
\label{proofharmonic1}
f(O)&\leq {c}_k\sum\limits_{i:J_i\subseteq S^l\setminus O}\varrho_{J_i}(S^{i-1})+\sum\limits_{i:J_i\subseteq O\cap S^l}\varrho_{J_i}(S^{i-1})+\sum\limits_{i:P_i\subseteq O\setminus S^l}\varrho_{P_i}(S^l)\nonumber\\
&={c}_k\sum\limits_{i:J_i\subseteq S^l\setminus O}\varrho_{J_i}(S^{i-1})+\sum\limits_{i:J_i\subseteq O\cap S^l}\varrho_{J_i}(S^{i-1})+\sum_{\substack{ i:P_i\subseteq O\setminus S^l\\i\neq l+1}}\varrho_{P_i}(S^l)\nonumber\\
\quad\quad&\quad+\varrho_{P_{l+1}}(S^l).
\end{align}
By the monotoneity of the set function $f$ and because $P_{i_j}\subseteq J_{i_j}'$ for $1\leq j\leq r$, we have
\begin{equation}
\label{proofharmonic2}
\varrho_{P_{i_j}}(S^l)\leq \varrho_{J_{i_j}'}(S^l)
\end{equation}
for $1\leq j\leq r$.
Based on the fact that $J_i'\subseteq O$ for $1\leq i\leq l+1$, and because $P_{i_j}\subseteq J_{i_j}'$ and $P_{i_j}\subseteq O\setminus S^l$, we have 
\begin{equation}
\label{proofharmonic3}
J_{i_j}'\subseteq O\setminus S^l.
\end{equation}
Combining (\ref{proofharmonic1})-(\ref{proofharmonic3}) results in
 \begin{align}
\label{proofharmonic4}
f(O)\leq&{c}_k\sum\limits_{i:J_i\subseteq S^l\setminus O}\varrho_{J_i}(S^{i-1})+\sum\limits_{i:J_i\subseteq O\cap S^l}\varrho_{J_i}(S^{i-1})+\sum_{\substack{ i:J_i'\subseteq O\setminus S^l\\i\neq l+1}}\varrho_{J_i'}(S^l)\nonumber\\
\quad\quad&+\varrho_{J_{l+1}'}(S^l).
\end{align}
By submodularity (see inequality (\ref{eqn:submodularimplies})), we have 
\begin{equation}
\label{proofharmonic5}
\varrho_{J_i'}(S^l)\leq \varrho_{J_i'}(S^{i-1})
\end{equation}
for $1\leq i\leq l+1$.
By statement~a in Lemma~\ref{lemma3}, we have
\begin{equation}
\label{proofharmonic6}
\varrho_{J_i'}(S^{i-1})\leq \varrho_{J_i}(S^{i-1})
\end{equation}
 for $1\leq i\leq l+1$.
By combining inequalities~(\ref{proofharmonic5}) and (\ref{proofharmonic6}), we have
\begin{equation}
\label{proofharmonic7}
\varrho_{J_i'}(S^l)\leq \varrho_{J_i}(S^{i-1}),
\end{equation}
for $1\leq i\leq l+1$. Combining inequalities~(\ref{proofharmonic4}) and (\ref{proofharmonic7}) results in 
 \begin{align}
\label{proofharmonic8}
f(O)&\leq  {c}_k\sum\limits_{i:J_i\subseteq S^l\setminus O}\varrho_{J_i}(S^{i-1})+\sum\limits_{i:J_i\subseteq O\cap S^l}\varrho_{J_i}(S^{i-1})+\sum_{\substack{ i:J_i'\subseteq O\setminus S^l\\i\neq l+1}}\varrho_{J_i}(S^{i-1})\nonumber\\
\quad\quad\quad&+\varrho_{J_{l+1}}(S^l).
\end{align}
By statement~b in Lemma~\ref{lemma3}, and inequality~(\ref{proofharmonic8}), we have 
 \begin{align}
\label{proofharmonic9}
f(O)&\leq {c}_k\sum\limits_{i:J_i\subseteq S^l\setminus O}\varrho_{J_i}(S^{i-1})+\sum\limits_{i:J_i\subseteq O\cap S^l}\varrho_{J_i}(S^{i-1})+\sum_{ i:J_i\subseteq S^l\setminus O}\varrho_{J_i}(S^{i-1})\nonumber\\
\quad\quad\quad&+\varrho_{J_{l+1}}(S^l).
\end{align}
Because $$\sum\limits_{i:J_i\subseteq S^l\setminus O}\varrho_{J_i}(S^{i-1})\leq f(S^l),$$ 

$$\sum\limits_{i:J_i\subseteq O\cap S^l}\varrho_{J_i}(S^{i-1})+\sum_{ i:J_i\subseteq S^l\setminus O}\varrho_{J_i}(S^{i-1})=f(S^l),$$
and
$$\varrho_{J_{l+1}}(S^l)=f(S)-f(S^l),$$
we can use inequality~(\ref{proofharmonic9}) to write
 \begin{align*}
f(O)\leq{c}_k f(S^l)+f(S^l)+f(S)-f(S^l)\leq ({c}_k+1)f(S),
\end{align*}
which implies that ${f(S)}/{f(O)}\geq \frac{1}{1+c_k}.$ \qed
\begin{remark}
For $k=1$, the harmonic bound for a general matroid becomes the bound in Theorem~\ref{Theorem2.1}.
\end{remark}
\begin{remark}
 The function $g(x)=1/(1+x)$ is nonincreasing in $x$ on the interval $]0,1]$. 
\end{remark}
\begin{remark}
The harmonic bound $1/(1+c_k)$ for the $k$-batch greedy strategy holds for \emph{any} matroid.  For the special case of a uniform matroid, we will give a different (exponential) bound in Theorem~\ref{Theorem3.4}  below. We will also show that this exponential bound is better than the harmonic bound when $k$ divides the rank of the uniform matroid $K$.
\end{remark}

In  Theorem~\ref{Theorem3.4} below, we provide an exponential bound for the $k$-batch greedy strategy in the case of uniform matroids. The special case when $c_k=1$ was derived in \cite{nemhauser19781}. Our result here is more general, and the method used in our proof is different from that of  \cite{nemhauser19781}. The new proof here is of particular interest because the technique here is not akin to that used in the case of general matroids in Theorem~\ref{Theorem3.3} and also  was not considered in \cite{nemhauser19781}. Before stating the theorem, we first present a proposition that  will be used in proving  Theorem~\ref{Theorem3.4}.

Choose a set $J^*\subseteq X\setminus S^l$ with $|J^*|=k$ so as to maximize  $f(S^l\cup J^*)-f(S^l).$ Write $\varrho_{J_{l+1}}(S^l)=f(S^{l+1})-f(S^l)\ \text{and}\ \varrho_{J^*}(S^l)=f(S^l\cup J^*)-f(S^l).$ We have the following proposition.
\begin{proposition}
\label{Prop2}
 Let $f:2^X\longrightarrow {\rm I\!R}$ be a submodular set function. Then when $(X,\mathcal{I})$ is a uniform matroid, we have $\varrho_{J_{l+1}}(S^l)\geq \frac{m}{k}\varrho_{J^*}(S^l)$.
 \end{proposition}
{\it Proof}  
  Let $\{M_1,\ldots, M_r\}$, where
\[r=\binom{k}{m},\]  be the collection of all the subsets  of $J^*$ with cardinality $m$. Then, each element of $J^*$ appears in exactly $p$ of these subsets, where 
 \[p=\binom{k-1}{m-1}.\] 
Using Lemma~\ref{lemma1} with $A=S^l, B=J^*$ and $B\setminus A=J^*$, we have 
\begin{equation}
\label{ineq1prop2}
 \sum\limits_{i=1}^r\varrho_{M_i}(S^l)\geq p\varrho_{J^*}(S^l).
\end{equation}
Because $|S^l\cup M_i|=kl+m= K$, by the definition of the uniform matroid $(X,\mathcal{I})$, we have $S^l\cup M_i\in\mathcal{I}.$
 By the definition of the $k$-batch greedy strategy and the monotoneity  of the set function $f$, we have
\begin{equation}
\label{ineq2prop2}
 \varrho_{J_{l+1}}(S^l)\geq \varrho_{M_i}(S^l),
\end{equation}
which implies that 
\begin{equation}
\label{ineq3prop2}
r\varrho_{J_{l+1}}(S^l)\geq \sum\limits_{i=1}^r\varrho_{M_i}(S^l).
\end{equation}
 Combining (\ref{ineq3prop2}) and (\ref{ineq1prop2}), we have 
 \[\varrho_{J_{l+1}}(S^l)\geq\frac{1}{r}\sum\limits_{i=1}^r\varrho_{M_i}(S^l)\geq \frac{p}{r}\varrho_{J^*}(S^l)\geq \frac{m}{k}\varrho_{J^*}(S^l),\] which implies that  $\varrho_{J_{l+1}}(S^l)\geq \frac{m}{k}\varrho_{J^*}(S^l)$.
 \qed

\begin{remark}
The reason we require $(X,\mathcal{I})$ to be a uniform matroid is that this result does not necessarily hold for a general matroid, because $S^l\cup M_i\in \mathcal{I}$ is not guaranteed for a general matroid, and in consequence inequality~(\ref{ineq2prop2}) does not necessarily hold.
\end{remark}
\begin{theorem}
\label{Theorem3.4}
Let $(X,\mathcal{I})$ be a uniform matroid and $f:2^X\longrightarrow {\rm I\!R}$ be a polymatroid set function. Then, any $k$-batch greedy solution $S$ satisfies 
\begin{align}
\label{k-batchuniformbound}
\frac{f(S)}{f(O)}&\geq \frac{1}{c_k}\left(1-\left(1-\frac{c_k}{l+1}\frac{m}{k}\right)\left(1-\frac{c_k}{l+1}\right)^l\right). 
\end{align}
\end{theorem}
\vspace{2mm}
{\it Proof}
Recall again that $J_i$ is the set selected at stage $i$ by the $k$-batch greedy strategy, $S^i=\cup_{j=1}^i J_j$ for $1\leq i\leq l$, and $S^0=\emptyset$ as defined in Section~3.1. Also recall that we defined $J^*$ as the set that maximizes  $f(S^l\cup J^*)-f(S^l)$ with $J^*\subseteq X\setminus S^l$ and $|J^*|=k$. 

Let $\{P_{i,1},\ldots, P_{i,r_i}\}$ be a partition of $O\setminus S^i$ satisfying $P_{i,j}\subseteq J_{i,j}'$ for $1~\leq~ j\leq r_i$. Finding $\{P_{i,1},\ldots, P_{i,r_i}\}$ for each $i$ is similar to finding $\{P_{i_1},\ldots, P_{i_r}\}$ which was given in the proof of Theorem~\ref{Theorem3.3}.
Letting $B=O$ and $A=S^i$ ($0\leq i\leq l$)  in Proposition~\ref{Prop1}, we have
\begin{equation}
\label{ineqexpproof1}
f(O\cup S^i)\leq f(S^i)+\sum\limits_{j: P_{i,j}\subseteq O\setminus S^i}\varrho_{P_{i,j}}(S^i).
\end{equation}
By the monotoneity of the set function $f$ and because  $P_{i,j}\subseteq J_{i,j}'$ for $1\leq j\leq r_i$, we have 
\begin{equation}
\label{ineqexpproof2}
\varrho_{P_{i,j}}(S^i)\leq \varrho_{J_{i,j}'}(S^i).
\end{equation}
Based on the fact that $J_{i,j}'\subseteq O$ and because $P_{i,j}\subseteq O\setminus S^i$ and $P_{i,j}\subseteq J_{i,j}'$, we have 
\begin{equation}
\label{ineqexpproof3}
J_{i,j}'\subseteq O\setminus S^i.
\end{equation}
Combining (\ref{ineqexpproof1})-(\ref{ineqexpproof3}) results in
\begin{equation}
\label{ineqexpproof4}
f(O\cup S^i)\leq f(S^i)+\sum\limits_{j: J_{i,j}'\subseteq O\setminus S^i}\varrho_{J_{i,j}'}(S^i).
\end{equation}
For $0\leq i\leq l-1$, we have $|S^i\cup J_{i,j}'|\leq K$, which implies that $S^i\cup J_{i,j}'\in \mathcal{I}$ always holds. 
So for any given $i$ ($0\leq i\leq l-1$), by the definition of the $k$-batch  greedy strategy, we have 
\begin{align}
\label{ineqexpproof5}
\varrho_{J_{i,j}'}(S^i)\leq \varrho_{J_{i+1}}(S^i)
\end{align}
for any $J_{i,j}'\subseteq O\setminus S^i$.
Now consider $i=l$. For any $J_{l,j}'\subseteq O\setminus S^l$ with $|J_{l,j}'|=k$, by the definition of $J^*$ before Proposition~\ref{Prop2}, we have 
\[\varrho_{J_{l,j}'}(S^{l})\leq \varrho_{J^*}(S^l).\]
By the definition of $J^*$ and the monotoneity of the set function $f$, we have 
\[\varrho_{J_{l+1}'}(S^l)\leq \varrho_{J^*}(S^l).\]
Combining the two inequalities above, we have for any $J_{l,j}'\subseteq O\setminus S^l$,
\begin{equation}
\label{ineqexpproof6}
\varrho_{J_{l,j}'}(S^l)\leq \varrho_{J^*}(S^l).
\end{equation}
By inequalities~(\ref{ineqexpproof5}) and (\ref{ineqexpproof6}), for any given $i$ ($0\leq i\leq l$), we have 
\begin{equation}
\label{ineqexpproof7}
\varrho_{J_{i,j}'}(S^i)\leq \varrho_{L_{i+1}}(S^i)
\end{equation}
for any $J_{i,j}'\subseteq O\setminus S^i$,
where 
\[ L_{i+1}=
    \begin{cases}
    J_{i+1}, & 0\leq i\leq l-1,\\
     J^*, & i=l.
    \end{cases}\]
By inequalities~(\ref{ineqexpproof4}) and~(\ref{ineqexpproof7}), we have 
\[f(O\cup S^i)\leq f(S^i)+\sum\limits_{j: J_{i,j}'\subseteq O\setminus S^i}\varrho_{L_{i+1}}(S^i),\]
which implies that 
\begin{equation}
\label{ineqexpproof8}
f(O\cup S^i)\leq f(S^i)+(l+1)\varrho_{L_{i+1}}(S^i).
\end{equation}
Setting $i=0$ in inequality~(\ref{ineqexpproof8}), recalling that $S^0=\emptyset$, and because $S^1=J_1$ by definition, we have 
\begin{equation}
\label{ineqexpproof9}
f(S^1)\geq \frac{1}{l+1}f(O).
\end{equation}
For $1\leq i\leq l$, we write 
\begin{equation}
\label{ineqexpproof10}
\frac{f(O\cup S^i)-f(O)}{f(S^i)}=\frac{\sum\limits_{j=0}^{i-1}f(O\cup S^j\cup J_{j+1})-f(O\cup S^j)}{\sum\limits_{j=0}^{i-1}f(S^j\cup J_{j+1})-f(S^j)}.
\end{equation}
By submodularity (see~(\ref{eqn:submodularimplies})), we have
\begin{equation}
\label{ineqexpproof11}
f(O\cup S^j\cup J_{j+1})-f(O\cup S^j)\geq f(X)-f(X\setminus J_{j+1})
\end{equation}
and
\begin{equation}
\label{ineqexpproof12}
f(S^j\cup J_{j+1})-f(S^j)\leq f(J_{j+1})-f(\emptyset)
\end{equation}
 for $0\leq j\leq i-1$.
By the definition of the total $k$-batch curvature, we have  
\begin{equation}
\label{ineqexpproof13}
\frac{f(X)-f(X\setminus J_{j+1})}{f(J_{j+1})-f(\emptyset)}\geq 1-c_k
\end{equation}
for $0\leq j\leq i-1$.
Combining inequalities~(\ref{ineqexpproof10})-(\ref{ineqexpproof12}) results in
\[\frac{f(O\cup S^i)-f(O)}{f(S^i)}\geq \frac{\sum\limits_{j=0}^{i-1}f(X)-f(X\setminus J_{j+1} )}{\sum\limits_{j=0}^{i-1}f( J_{j+1})-f(\emptyset)}\geq 1-c_k.\]
This in turn implies that
\[f(O)+(1-{c}_k)f(S^i)\leq f(O\cup S^i).\] 
Combining the above  inequality and (\ref{ineqexpproof8}), we have 
\begin{equation}
\label{ineqexpproof14}
f(S^i\cup L_{i+1})\geq \frac{1}{l+1}f(O)+\left(1-\frac{{c}_k}{l+1}\right)f(S^i)
\end{equation}
for $1\leq i\leq l$.
By inequality~(\ref{ineqexpproof9}) and successive application of inequality~(\ref{ineqexpproof14}) for $i=1,\ldots, l$, we have
\begin{align}
\label{ineq:uniform1}
f(S^{l})&\geq \frac{1}{l+1}f(O)+\left(1-\frac{{c}_k}{l+1}\right)f(S^{l-1})\nonumber\\
&\geq \frac{1}{l+1}f(O)\sum\limits_{i=0}^{l-1}\left(1-\frac{{c}_k}{l+1}\right)^i\nonumber\\
&= \frac{1}{{c}_k}\left(1-\left(1-\frac{{c}_k}{l+1}\right)^{l}\right)f(O),
\end{align}
and
\begin{align}
\label{ineq:uniform2}
f(S^{l}\cup J^*)&\geq \frac{1}{l+1}f(O)+\left(1-\frac{c_k}{l+1}\right)f(S^l)\nonumber\\
&\geq \frac{1}{l+1}f(O)\sum\limits_{i=0}^l\left(1-\frac{c_k}{l+1}\right)^i\nonumber\\
&\geq \frac{1}{{c}_k}\left(1-\left(1-\frac{{c}_k}{l+1}\right)^{l+1}\right)f(O).
\end{align}
Using Proposition~\ref{Prop2} and combining inequalities~(\ref{ineq:uniform1}) and~(\ref{ineq:uniform2}), we have
\begin{align*}
f(S) &\geq  \frac{m}{k}f(S^l\cup J^*)+\left(1-\frac{m}{k}\right)f(S^l)\\
&\geq\frac{m}{k}\frac{1}{{c}_k}\left(1-\left(1-\frac{{c}_k}{l+1}\right)^{l+1}\right)f(O)+\\
&\quad \quad\quad \left(1-\frac{m}{k}\right)\frac{1}{{c}_k}\left(1-\left(1-\frac{{c}_k}{l+1}\right)^{l}\right)f(O)\\
&= \frac{1}{{c}_k}\left(1-\left(1-\frac{{c}_k}{l+1}\frac{m}{k}\right)\left(1-\frac{c_k}{l+1}\right)^l\right)f(O),
\end{align*}
which implies (\ref{k-batchuniformbound}).
\qed
\begin{remark}
 For $k=1$, the exponential bound for a uniform matroid becomes the bound in Theorem~\ref{Theorem2.2}.
\end{remark}
\begin{remark}
 The exponential bound for a uniform matroid becomes $$1-(1-\frac{1}{l+1}\frac{m}{k})(1-\frac{1}{l+1})^l$$ for $c_k=1$,  which is the bound  in \cite{nemhauser19781}.
\end{remark}
\begin{remark}
When $m=k$, i.e., when $k$ divides the cardinality $K$,  the exponential bound for a uniform matroid becomes $\frac{1}{c_k}(1-(1-\frac{c_k}{l+1})^{l+1})$, which is the bound in \cite{Liu2016}. 
\end{remark}
\begin{remark}
 Let $g(x,y)= \frac{1}{x}\left(1-(1-\frac{x}{y})^{y}\right)$. The function $g(x, y)$ is nonincreasing in $x$ on the interval $]0, 1]$ for
any positive integer $y$. Also, $g(x, y)$ is nonincreasing in $y$ when $x$ is a constant on the interval $]0, 1]$.
\end{remark}
\begin{remark}
 Even if the total curvature $c_k$ is monotone in $k$, the exponential bound for a uniform matroid is not necessarily monotone. But under the condition that $k$ divides $K$, it is monotone. To be specific, if $k$ divides $K$, then $K=k(l+1)$ for some positive integer $l$. Thus, as $k$ increases, $l+1$ decreases, and if $c_k$ decreases, then, we have that $\frac{1}{c_k}(1-(1-\frac{c_k}{l+1})^{l+1})$ is nondecreasing in $k$ based on the previous remark.
\end{remark}
\begin{remark}
 {When $m=k$, the exponential bound is tight, as shown in \cite{nemhauser19781}. Moreover, for this case, the exponential bound $\frac{1}{{c}_k}(1-(1-\frac{{c}_k}{l+1}\frac{m}{k})(1-\frac{c_k}{l+1})^l)$ is better than the harmonic bound $1/(1+c_k)$ because $$(1-(1-\frac{c_k}{l+1})^{l+1})/c_k\geq (1-e^{-c_k})/c_k\geq 1/(1+c_k).$$ However, if $k$ does not divide $K$, the exponential bound might be worse than the harmonic bound. For example, when $K=100, k=80$, and $c_k=0.6$,  the exponential bound is $0.5875$, which is worse than the harmonic bound  $0.6250$.}
\end{remark}
\begin{remark}
 The monotoneity of $1/(1+c_k)$  implies that the $k$-batch greedy strategy has a better harmonic bound  than the $1$-batch greedy strategy if $c_k\leq c$. The monotoneity of  $\frac{1}{c_k}(1-(1-\frac{c_k}{l+1})^{l+1})$ implies that the $k$-batch ($k$ divides $K$) greedy strategy has a better exponential bound than the $1$-batch  greedy strategy if $c_k\leq c$.
\end{remark}

The following theorem establishes that indeed $c_k\leq c$.

\begin{theorem}
\label{Theorem3.5}
Let $f:2^X\longrightarrow {\rm I\!R}$ be a polymatroid set function with total curvature $c$ and total $k$-batch curvatures $\{c_k\}_{k=1}^K$. Then, $c_k\leq c$ for $1\leq k\leq K$.
\end{theorem}
{\it Proof}
By the definition of the total $k$-batch curvature $c_k$, we have 
\begin{align*}
c_k=\max_{I\in \hat{X}}\left\{1-\frac{\varrho_{J}({X\setminus I})}{\varrho_{I}(\emptyset)}\right\}=1-\min_{I\in\hat{X}}\left\{\frac{\sum\limits_{j=1}^k\varrho_{i_j}(X\setminus I_j)}{\sum\limits_{j=1}^k\varrho_{i_j}(I_{j-1})}\right\},
\end{align*}
where  $I=\{i_1,\ldots, i_k\}$ and $I_j=\{i_1,\ldots, i_j\}$ for $1\leq j\leq k$.

By submodularity (see~(\ref{eqn:submodularimplies})), we have
$$\varrho_{i_j}(X\setminus I_j)\geq \varrho_{i_j}(X\setminus\{i_j\})
\ \text{and}\ \varrho_{i_j}(I_{j-1})\leq \varrho_{i_j}(\emptyset)$$
 for $1\leq j\leq k$, which imply that 
\[\frac{\sum\limits_{j=1}^k\varrho_{i_j}(X\setminus I_j)}{\sum\limits_{j=1}^k\varrho_{i_j}(I_{j-1})}\geq\frac{\sum\limits_{j=1}^k\varrho_{i_j}(X\setminus\{i_j\})}{\sum\limits_{j=1}^k\varrho_{i_j}(\emptyset)}. \]
Therefore, we have
\begin{equation}
\label{Inequality2}
c_k\leq 1-\min_{I_k\in\hat{X}}\left\{\frac{\sum\limits_{j=1}^k\varrho_{i_j}(X\setminus\{i_j\})}{\sum\limits_{j=1}^k\varrho_{i_j}(\emptyset)}\right\}.
\end{equation}
By the definition of $c$ and the fact that $f$ is a polymatroid set function, we have 
$\varrho_{i_j}(X\setminus\{i_j\})\geq (1-c)\varrho_{i_j}(\emptyset)$
for $1\leq i\leq k$.
Combining this inequality and (\ref{Inequality2}), we have $c_k\leq 1-(1-c)=c.$
\qed
One would  expect the following generalization of Theorem~\ref{Theorem3.5} to hold: if $k_2\geq k_1$, then $c_{k_2}\leq c_{k_1}$. In the case of general matroid constraints, this conclusion implies that the bound is nondecreasing in $k$. In the case of uniform matroid constraints, monotoneity of the bound holds under the condition that $k$ divides $K$. We now state and prove the following theorem on the monotoneity of $c_k$, using Lemmas~\ref{lemma1} and~\ref{lemma2} (Lemmas~1.1 and~1.2 in \cite{vondrak2010}).

\begin{theorem}
\label{Theorem3.6}
Let $f:2^X\longrightarrow {\rm I\!R}$ be a polymatroid set function with $k$-batch curvatures $\{c_k\}_{k=1}^K$. Then, $c_{k_2}\leq c_{k_1}$ whenever $k_2\geq k_1$.
\end{theorem}
{\it Proof} 
Let $J\subseteq X$ be a set with cardinality $k_2$ satisfying $f(J)>0$. Let $\{M_1,\ldots, M_s\}$ be the collection of all the subsets  of $J$ with cardinality $k_1$ ($k_1\leq k_2$), where 
 \[s=\binom{k_2}{k_1}.\]
Then, 
each element of $J$ appears in exactly $q$ of the subsets $\{M_1,\ldots, M_s\}$, where
\[q=\binom{k_2-1}{k_1-1}.\]
Using Lemma~\ref{lemma2} with $A=X$,  $A'=X\setminus J$, and $A\setminus A'=J$, we  have 
\[\sum\limits_{i=1}^s\varrho_{M_i}(X\setminus M_i)\leq q \varrho_J(X\setminus J),\] which implies that
\begin{equation}
\label{comparison1}
 \varrho_J(X\setminus J)\geq \frac{1}{q} \sum\limits_{i=1}^s\varrho_{M_i}(X\setminus M_i).
\end{equation}
Based on the fact that $\{M_1,\ldots, M_s\}$ is the collection of all the subsets of $J$  with cardinality $k_1$ and that
each element of $J$ appears in exactly $q$ of these subsets, using Lemma~\ref{lemma1} with $B=J$ and $A=\emptyset$, we have
\[\sum\limits_{i=1}^s\varrho_{M_i}(\emptyset)\geq q \varrho_J(\emptyset),\] which implies that
\begin{equation}
\label{comparison2}
\varrho_J(\emptyset)\leq \frac{1}{q}\sum\limits_{i=1}^s\varrho_{M_i}(\emptyset).
\end{equation}
Combining inequalities~(\ref{comparison1}) and~(\ref{comparison2}) results in 
\begin{align}
\label{fraction1}
\frac{ \varrho_J(X\setminus J)}{\varrho_J(\emptyset)}&\geq \frac{ \frac{1}{q} \sum\limits_{i=1}^s\varrho_{M_i}(X\setminus M_i)}{\frac{1}{q}\sum\limits_{i=1}^s\varrho_{M_i}(\emptyset)}=\frac{  \sum\limits_{i=1}^s\varrho_{M_i}(X\setminus M_i)}{\sum\limits_{i=1}^s\varrho_{M_i}(\emptyset)}.
\end{align}
Recall the definition of the total $k$-batch curvature $c_{k}$ in~(\ref{totalkbatchcurvature}). Because $|M_i|=k_1$ for $1\leq i\leq s$ and $f$ is a polymatroid set function, we have 
\begin{equation}
\label{curvatureineq}
\varrho_{M_i}(X\setminus M_i)\geq (1-c_{k_1})\varrho_{M_i}(\emptyset)
\end{equation}
 for $1\leq i\leq s$.
Combining inequalities~(\ref{fraction1}) and~(\ref{curvatureineq}) results in 
\begin{equation}
\label{lastinequal}
\frac{ \varrho_J(X\setminus J)}{\varrho_J(\emptyset)}\geq\frac{ (1-c_{k_1})\sum\limits_{i=1}^s\varrho_{M_i}(\emptyset)}{\sum\limits_{i=1}^s\varrho_{M_i}(\emptyset)}=1-c_{k_1}.
\end{equation}
By~(\ref{totalkbatchcurvature}), $c_{k_2}$ can be written as 
\begin{equation}
\label{ck2curvature}
c_{k_2}=1-\min_{J\in\hat{X}}\left\{\frac{\varrho_J(X\setminus J)}{\varrho_J(\emptyset)}\right\}.
\end{equation}
By (\ref{lastinequal}) and (\ref{ck2curvature}), we have $c_{k_2}\leq 1-(1-c_{k_1})=c_{k_1}.$
\qed
\begin{remark}
When $k_1=1$ and $k_2=k$, Theorem~\ref{Theorem3.6} reduces to Theorem~\ref{Theorem3.5}. However, the proof of Theorem~\ref{Theorem3.5} can be  used only to prove the case when $k_1$ divides $k_2$ in Theorem~\ref{Theorem3.6}. This is why we have chosen to separate the two theorems.
\end{remark}
\section{Examples}\label{sc:IV}
In this section, we consider a task scheduling problem and an adaptive sensing problem to illustrate our results. Specially, we demonstrate that the total curvature $c_k$ decreases in $k$ and the performance bound for a uniform matroid  increases in $k$ under the condition that $k$ divides $K$.

\subsection{Task Scheduling}
As a canonical example for problem (\ref{eqn:1}), we  consider the task scheduling problem  posed in \cite{streeter2008online}, which was also analyzed in \cite{ZhC13J} and \cite{YJ2015}.  In this problem, there are $n$ subtasks and a  set $X$ of $N$ agents. At each stage, a subtask $i$ is assigned to an agent $a$, who accomplishes the task with probability $p_i(a)$. Let $X_i(\{a_1,a_2,\ldots, a_k\})$ denote the Bernoulli random variable that signifies whether or not subtask $i$ has been accomplished after performing the set of agents $\{a_1,a_2,\ldots, a_k\}$ over $k$ stages. Then $\frac{1}{n}\sum_{i=1}^n	X_i(\{a_1,a_2,\ldots,a_k\})$ is the fraction of subtasks accomplished after $k$ stages by employing agents $\{a_1,a_2,\ldots, a_k\}$. The objective function $f$ for this problem is the expected value of this fraction, which can be written as
$$f(\{a_1,\ldots,a_k\})=\frac{1}{n}\sum_{i=1}^n\left(1-\prod_{j=1}^k\left(1-p_i(a_j)\right)\right).$$
Assume that $p_i(a)>0$ for any $a\in X$. Then it is easy to check that $f$ is nondecreasing. Therefore, when $\mathcal{I}=\{S\subseteq X: |S|\leq K\}$, this problem has an optimal solution of length $K$.  Also, it is easy to check that  $f$ has the diminishing-return property and $f(\emptyset)=0$. Thus, $f$ is a polymatroid set function.

For convenience, we only consider the special case $n=1$; our analysis can be generalized to any $n\geq 2$. For $n=1$, we have 
$$f(\{a_1,\ldots,a_k\})=1-\prod_{j=1}^k\left(1-p(a_j)\right),$$
where $p(\cdot)=p_1(\cdot)$.

Let us  order the elements of $X$ as $a_{[1]}, a_{[2]}, \ldots, a_{[N]}$ such that $$0<p(a_{[1]})\leq p(a_{[2]})\leq\ldots\leq p(a_{[N]})\leq 1.$$ Then by the definition of the total curvature $c_k$, we have 
\[
c_k=\max\limits_{i_1,\ldots,i_k\in {X}}\left\{1-\frac{f(X)-f(X\setminus\{i_1,\ldots,i_k\})}{f(\{i_1,\ldots, i_k\})-f(\emptyset)}\right\}=1-\prod_{l=k+1}^N(1-p(a_{[l]})).
\]
%\begin{figure}
%   \centering
%   \includegraphics[width=\columnwidth, trim=50 129 20 130, clip]{FIG1.pdf}
%\label{Fig1}
%\caption{Task scheduling example}
%\end{figure}
To numerically evaluate the relevant quantities here, we randomly generate a set of $\{p(a_i)\}_{i=1}^{30}$.
 In Figure~1, we consider $ K=20$,  and  batch sizes $k=1,2,\ldots, 10$.  From the expression of $c_k$, we can see that $c_k$ is nonincreasing in $k$, but when $N$ is large, $c_k$ is close to 1 for each $k$. Figure~1 shows that  the exponential bound for $k=3,6,8,9$ is worse than that for $k=1,2$, which illustrates our earlier remark that the exponential bound for the uniform matroid case is not necessarily 
e in $k$  even though $c_k$ is monotone in $k$. Figure~1 also shows that  the exponential bound $\frac{1}{{c}_k}(1-(1-\frac{{c}_k}{l+1}\frac{m}{k})(1-\frac{c_k}{l+1})^l$ coincides with $\frac{1}{c_k}(1-(1-\frac{c_k}{l+1})^{l+1})$ for $k=1,2,4,5,10$ and it is nondecreasing in $k$, which illustrates our remark that  the exponential bound is nondecreasing in $k$ under the condition that $k$ divides $K$.
\begin{figure}[ht]
\begin{center}
\vspace{-1in}
\includegraphics[width=3.0 in]{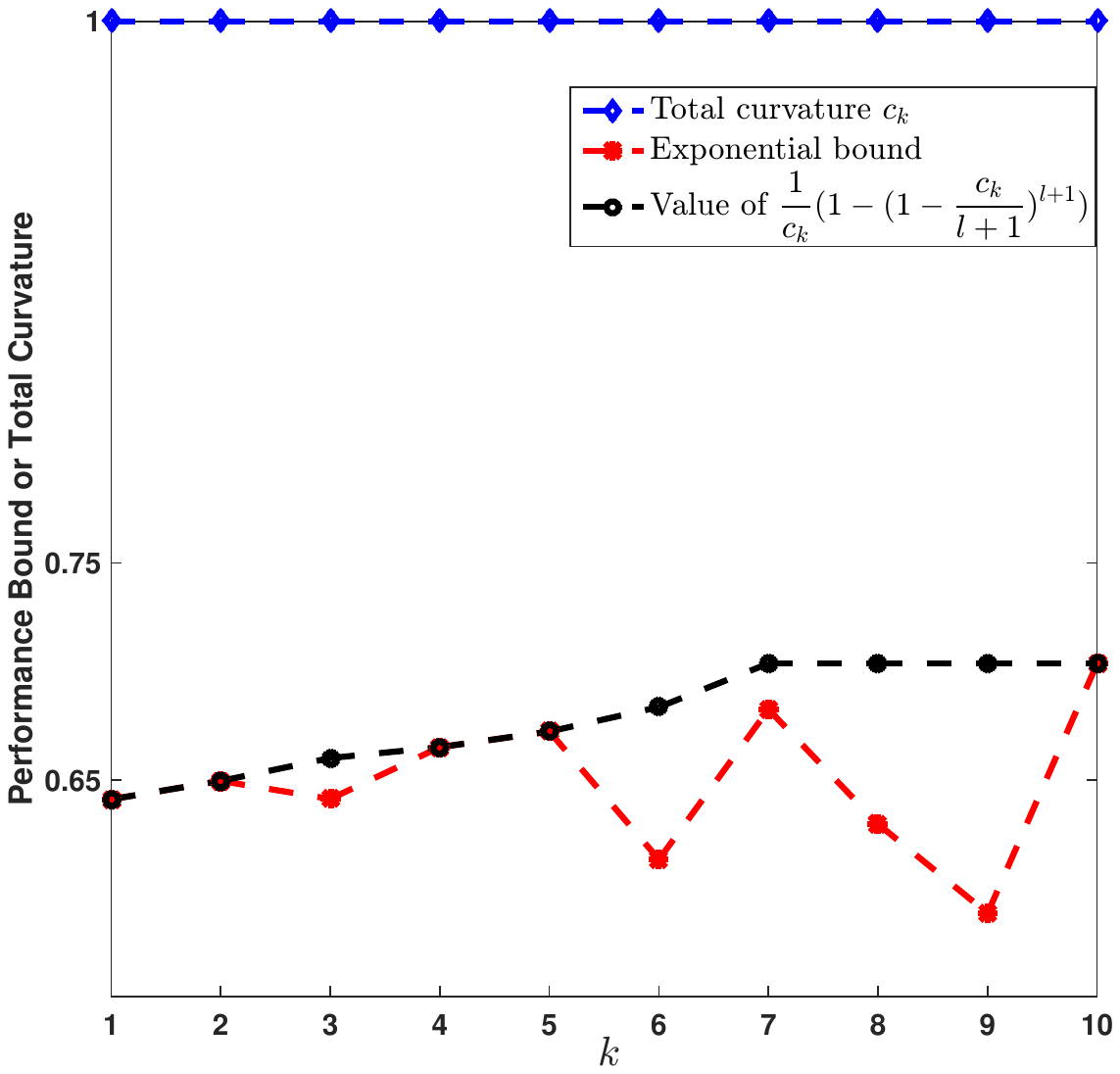}
% The figure caption is below the figure.
\end{center}
\vspace{-0.8in}
\caption{Task scheduling example }
\end{figure}

Owing to the nature of the total curvature for this example, it is not easy to see that $c_k$ is nonincreasing in $k$ (all $c_k$ values here are very close to 1). The next example will illustrate that the total curvature does decrease in $k$ and again demonstrate our claim that the exponential bound for the uniform matroid case is not necessarily monotone in $k$ but it is monotone in $k$ under the condition that $k$ divides $K$.
\subsection{Adaptive Sensing}

As our second example application, we consider the adaptive sensing design problem posed in \cite{LiC12} and \cite{ZhC13J}. Consider a signal of interest $x \in {\rm I\!R}^2$ with normal prior distribution $\mathcal{N} (0, I)$, where $I$ is the $2\times 2$ identity matrix; our analysis easily generalizes to dimensions larger than $2$. Let $\mathbb{B}=\{\mathrm{Diag}(\sqrt{b},\sqrt{1-b}):  b\in\{b_1,\ldots,b_N\}\}$, where $ b_i\in[0.5,1]$ for $1\leq i
\leq N$. At each stage $i$, we make a measurement $y_i$ of the form
\[
y_i=B_ix+w_i, 
\]
where $B_i \in\mathbb{B}$ and $w_i$ represents i.i.d.\ Gaussian measurement noise with mean zero and covariance $\sigma^2 I$, independent of $x$.

The objective function $f$ for this problem is the information gain, which can be written as
\[
f(\{B_1,\ldots, B_k\})=H_0-H_k.
\]
Here, $H_0=\frac{N}{2}\text{log}(2\pi e)$ is the entropy of the prior
distribution of $x$ and $H_k$ is the entropy of the posterior
distribution of $x$ given $\{y_i\}_{i=1}^k$; that is,
\[
H_k=\frac{1}{2}\text{log det}(P_k)+\frac{N}{2}\text{log}(2\pi e),
\] 
where 
\[
P_k=\left(P_{k-1}^{-1}+\frac{1}{\sigma^2}B_k^TB_k\right)^{-1}
\]
is the posterior covariance of $x$ given $\{y_i\}_{i=1}^k$ \cite{LiC12}.

The objective is to choose a set of measurement matrices $\{B_i^*\}_{i=1}^K$, $B_i^*\in\mathbb{B}$, to maximize the information gain $
f(\{B_1,\ldots, B_K\})=H_0-H_K$. It is easy to check that $f$ is nondecreasing, submodular, and $f(\emptyset)=0$; i.e., $f$ is a polymatroid set function.

For convenience, let $\sigma=1$.
Then, we have 
\begin{align*}
c_k&=\max\limits_{J_k\subseteq X, |J_k|=k}\left\{1-\frac{f(X)-f(X\setminus J_k)}{f(J_k)}\right\}\nonumber\\
&=\max\limits_{J_k\subseteq X, |J_k|=k}\left\{1- \frac{\log (st)-\log\left(s-\sum\limits_{i: e_i\in J_k}e_i\right)\left(t-\sum\limits_{i: e_i\in J_k}(1-e_i)\right)}{\log\left(1+\sum\limits_{i: e_i\in J_k}e_i\right)\left(1+\sum\limits_{i: e_i\in J_k}(1-e_i)\right)}\right\},
\end{align*}
where $X=\{B_1,\ldots, B_N\}$, $s=1+\sum_{i=1}^Ne_i$, and $t=1+\sum_{i=1}^N(1-e_i)$.

\begin{figure}[ht]
\begin{center}
\vspace{-0.6in}
\includegraphics[width=3.0 in]{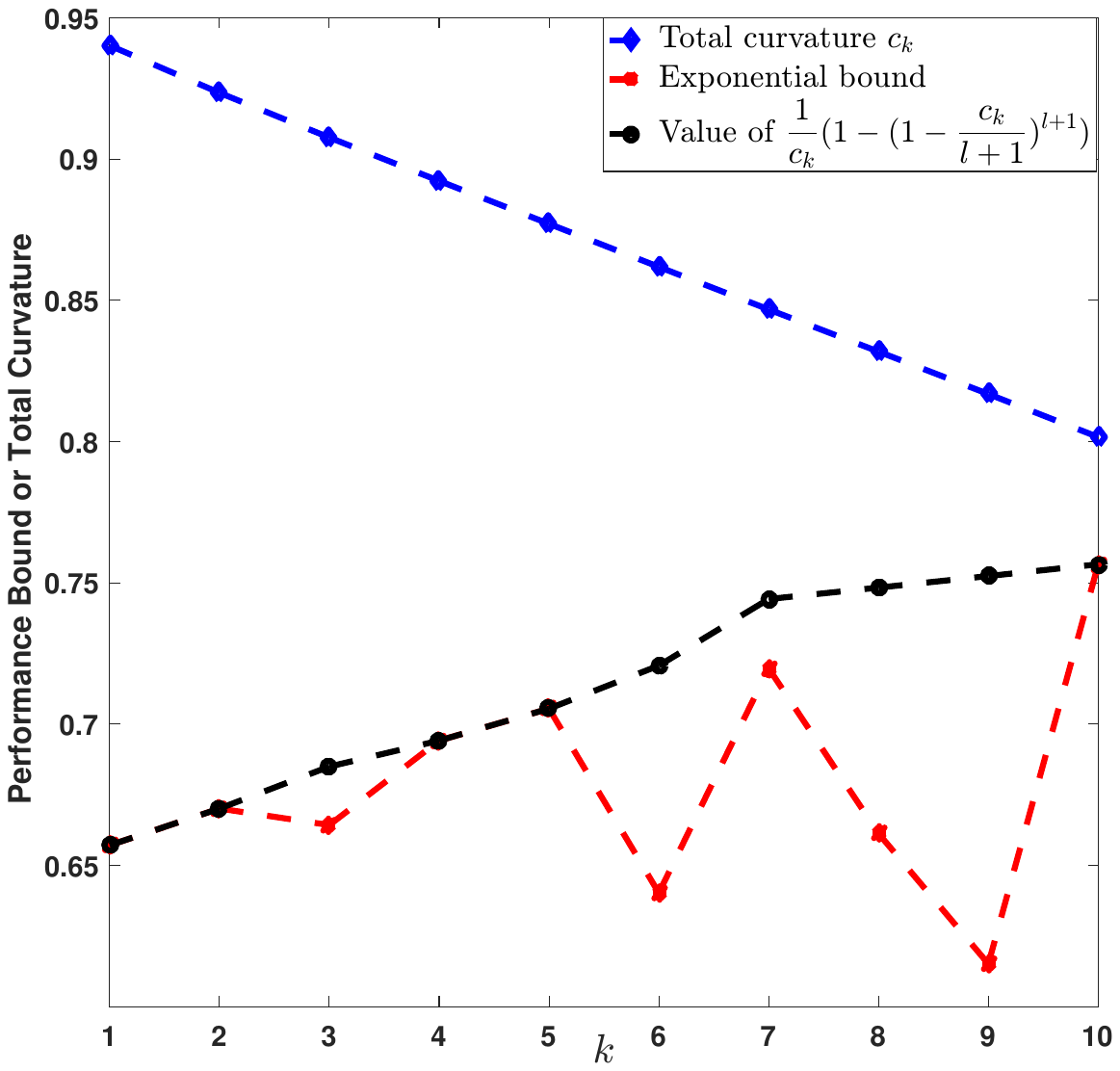}
% The figure caption is below the figure.
\end{center}
\vspace{-0.9in}
\caption{Adaptive sensing example}
\end{figure}

To numerically evaluate the relevant quantities here, we randomly generate a set of $\{e_i\}_{i=1}^{30}$. We first still consider $ K=20$ for $k=1,\ldots,10$ in Figure~2.
 Figure~2 shows that the total curvature decreases in $k$, while the exponential bound for the uniform matroid case only increases for $k=1, 2,4,5,7,10$ and the bound for $k=3, 6,8,9$ is worse than that for $k=1,2$. This illustrates that the exponential bound for the uniform matroid  case is not necessarily monotone in $k$.

Next, we consider $ K=24$ for $k=1,2,3,4,6,8$ in Figure~3. Figure~3 shows that the curvature decreases in $k$ and the exponential bound increases in $k$ since $k$ divides $K$ for $k=1,2,3,4,6,8$, which again demonstrates our claim that $c_k$ decreases in $k$ and the exponential bound increases in $k$ under the condition that $k$ divides $K$. \begin{figure}[ht]
\begin{center}
\vspace{-0.6in}
\includegraphics[width=3 in]{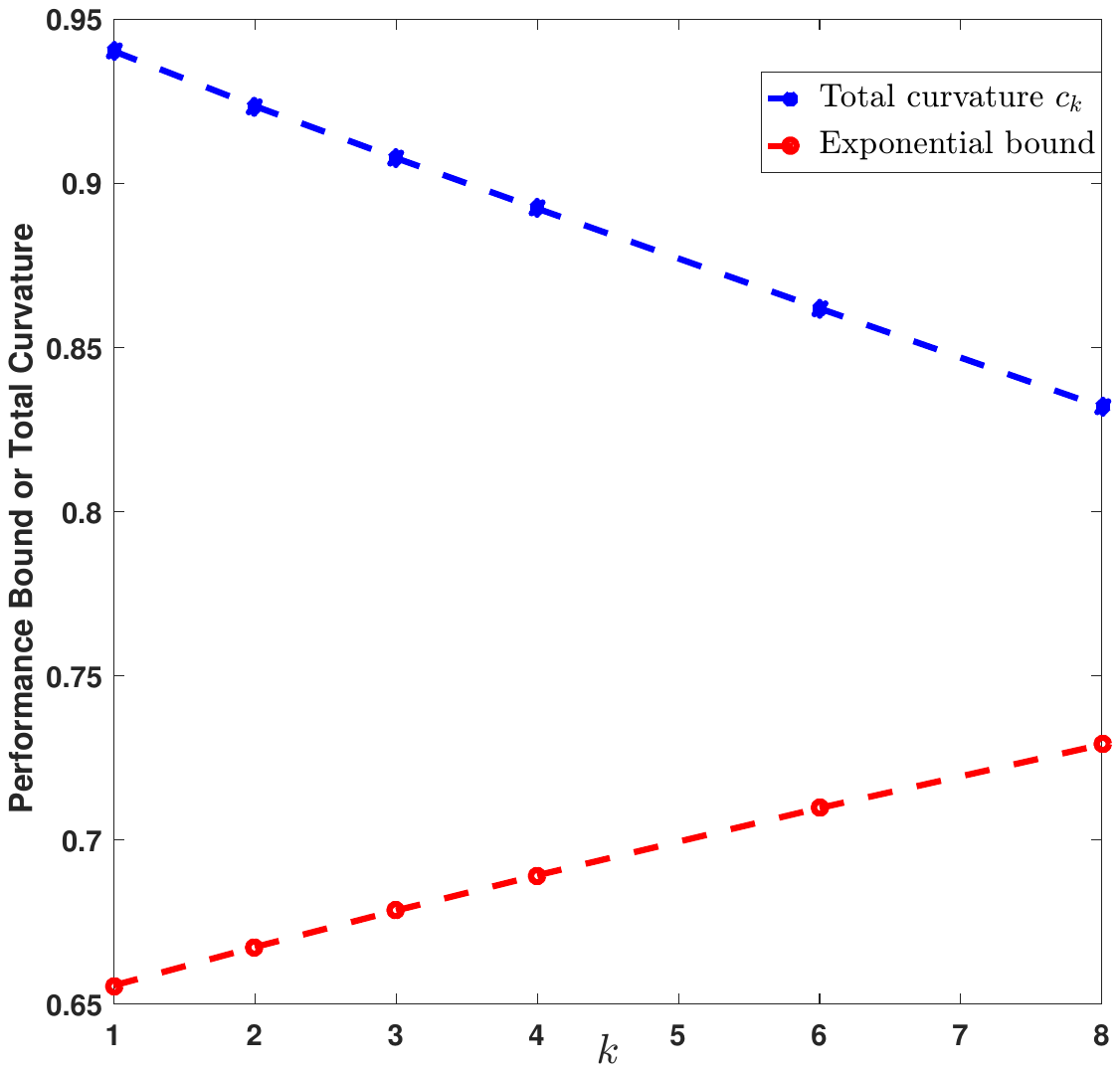}
% The figure caption is below the figure.
\end{center}
\vspace{-0.9in}
\caption{Adaptive sensing example }
\end{figure}

\section{Perspectives and Open Problems}
As defined, the total batched curvature depends on values of the objective function $f$ outside the constraint matroid $\mathcal{I}$. This raises two issues. First, the optimization problem~\eqref{eqn:1} makes sense even if the objective function $f$ is defined only on $\mathcal{I}$; the values of $f$ outside of $\mathcal{I}$ are irrelevant. Yet, if $f$ is only defined on $\mathcal{I}$, the total batched curvature is ill defined, which implies that the performance bounds we derived in this paper do not apply to the case when the objective function is defined only on $\mathcal{I}$. This suggests that a different definition of total batched curvature is possible. Second, even if $f$ were defined on the entire $2^X$, in which case the total batched curvature is well defined, it still remains that the values of $f$ outside of $\mathcal{I}$ are irrelevant to problem~\eqref{eqn:1}. This in turn suggests that they should also be irrelevant to the performance of greedy strategies. These reasons motivate a reformulation of the notion of total curvature that depends only on values of $f$ inside $\mathcal{I}$ and an investigation of the performance bounds in terms of the new curvature, and is part of our ongoing work.

A related problem setting in one where the argument of the objective function is not a set but an \emph{ordered} tuple, called a string.
The performance of  the $1$-batch greedy strategy for string optimization problem has been investigated in \cite{ZhC13J, streeter2008online}; however, the performance of the general $k$-batch greedy strategy for string optimization has not been investigated so far. As was the case in this paper, lifting does not work in the string setting. Moreover, batching in string submodular functions does not preserve submodularity in general. This makes analyzing the $k$-batch greedy strategy for string problems more challenging than for set problems, and remains open to date.

Another related problem arises in game theory. It turns out that similar techniques can be used to bound the performance of Nash equilibria in noncooperative games. The connection to the game setting is easy to imagine by associating our objective function with a social utility function, greedy strategies with Nash equilibria, and batching with cooperation of subgroups in games. In the game setting, the bounds are typically smaller (e.g., $1/2$ instead of $1-e^{-1}$), reflecting the \emph{price of anarchy}. However, similar to the result that batching often results in improved bounds, in the game setting it might be possible to show that cooperation of subgroups results in improved bounds. This is part of our ongoing work.

\section{Conclusions}
In this paper, we developed bounds on the performance of the batched greedy strategy relative to the optimal strategy in terms of a parameter called the total batched curvature. We showed that when the objective function is a polymatroid set function, the batched greedy strategy satisfies a harmonic bound for a general matroid constraint and an exponential bound for a uniform matroid constraint, both in terms of the total batched curvature.
We also studied the behavior of the bounds as functions of the batch size. Specifically, we proved that the harmonic bound for a general matroid is nondecreasing in the batch size and the exponential bound for a uniform matroid is nondecreasing in the batch size under the condition that the batch size divides the rank of the uniform matroid. Finally, we illustrated our results by considering a task scheduling problem and an adaptive sensing problem.
\begin{acknowledgements}
This work is supported in part by NSF under award CCF-1422658, and by the CSU Information Science and Technology Center (ISTeC). We would like to acknowledge the anonymous reviewers  for their insightful comments on our conference paper  \cite{Liu2016}, as these comments led us to improve our work.
\end{acknowledgements}

\appendix  %This command ends the counting of sections.
\section{Discussion on Matroid Preservation}
\label{AppendixA}
Suppose that $(X,\mathcal{I})$ is a uniform matroid.
In this appendix, we will provide an example to prove that the collection of subsets of $X$ of size $k$ satisfying the constraint $\mathcal{I}$ (i.e., actions in the $k$-batch greedy strategy) is \emph{not} in general a matroid.  This shows that lifting does not work; i.e., it is not in general possible to appeal to bounds for the $1$-batch greedy strategy to derive bounds for the $k$-batch greedy strategy. For convenience, we assume that $k$ divides the uniform matroid rank $K$.

Recall that for a matroid $(X,\mathcal{I})$, we have the following two properties:
\begin{itemize}
\item  [i.] For all $B\in\mathcal{I}$, any set $A\subseteq B$ is also in $\mathcal{I}$.
\item  [ii.] For any $A,B\in \mathcal{I}$, if $|B|>|A|$, then there exists $j\in B\setminus A$ such that $A\cup\{j\}\in\mathcal{I}$.
\end{itemize}

To apply lifting, first fix $k$. We will define a pair $(Y,\mathcal{J})$ such that $Y$ is the ``ground set''  of all $k$-element subsets of
$X$: $Y=\{y: y=\{a_1,\ldots, a_k\}, k $ is given, and $ a_i\in X\}$. Next, $\mathcal{J}$ is the set of all subsets of $Y$ such that their elements are disjoint and the union of their elements lies in $\mathcal{I}$.  
The following example shows that $(Y,\mathcal{J})$ constructed this way is not in general a matroid.

\begin{example}
\label{exampleA}
Fix $k=2$. Let $X=\{a,b,c,d\}$, and $\mathcal{I}$ be the power set of $X$ (a special case of a uniform matroid, with rank $K=4$). We have $Y=\{\{a,b\},\{a,c\},\{a,d\},$
$ \{b,c\},\{ b,d\},\{c,d\}\}$. 
Let $\mathcal{J}$ be as defined above.
\end{example}
We will now prove that $(Y,\mathcal{J})$ does not satisfy property~ii above. To see this, consider $A=\{\{a,b\}\}\in\mathcal{J}$ and $B=\{\{a,c\},\{b,d\}\}\in\mathcal{J}$. We have $|A|=1$ and $|B|=2$. 
Notice that $\{a,b\}\cap \{a,c\}\neq \emptyset$ and $\{a,b\}\cap \{b,d\}\neq \emptyset$. So, in this case clearly there does not exist $j\in B\setminus A$ such that $A\cup \{j\}\in\mathcal{J}$. Hence, property~ii fails and $(Y,\mathcal{J})$ is not a matroid.
\section{Comparing Different $k$-Batch Greedy Strategies}
\label{AppendixB}
It is tempting to think that the $k$-batch ($k\geq 2$) greedy strategy always outperforms the $1$-batch greedy strategy. In fact, this is false. To show this, we will provide two examples based on the \emph{maximum $K$-coverage problem}, which was considered in  \cite{Chandu2015}  to demonstrate via Monte Carlo simulations that the $1$-batch greedy strategy can perform better than the $2$-batch greedy strategy. The maximum $K$-coverage problem is to select at most $K$ sets from a collection of sets such that the  union of the selected sets has the maximum number of elements. Example~\ref{exampleB1} below is to choose at most 3 sets from a collection of 5 sets, and Example~\ref{exampleB2} is to choose at most 4 sets from a collection of 6 sets. In contrast to \cite{Chandu2015}, our examples are not based on simulation, but are analytical counterexamples.

\begin{example}  
\label{exampleB1}
Fix $K=3$ and let the sets to be selected be
$S_1=\{a,f\}$, $S_2=\{f\}$, $S_3~=~\{a,b,g\}$, $S_4=\{c,f,g\}$, and $S_5=\{e,g,h\}$.
\end{example}
For the $1$-batch greedy strategy, one solution is $\{S_3, S_4,S_5\}$, and the union of the selected sets is 
$S_3\cup S_4\cup S_5=\{a,b,c,e,f,g,h\}$.
For the $2$-batch greedy strategy, one solution is $\{S_1, S_5,S_3\}$, and the union of the selected sets is $S_1\cup S_5\cup S_3=\{a,b,e,f,g,h\}$.
It is easy to see that $|S_3\cup S_4\cup S_5|=7> |S_1\cup S_5\cup S_3|=6$.
\begin{example}
\label{exampleB2}
 Fix $K=4$ and let the sets to be selected be $S_1=\{h,i,j\},S_2=\{b,e,i,j\},$ $S_3=\{c,d,e,h\},S_4=\{b,d,f,h,i\}, S_5=\{a,h,i,j\}$, and $S_6=\{c,g,i\}$.
\end{example}
For the $1$-batch greedy strategy, one solution is $\{S_4, S_2,S_6, S_5\}$, and the union of the selected sets is 
$S_4\cup S_2\cup S_6\cup S_5=\{a,b,c,d,e,f,g,h,i,j\}$.
For the $2$-batch greedy strategy, one  solution is $\{S_2\cup S_3\cup S_4\cup S_5\}$, and the union of the selected sets is 
$S_2\cup S_3\cup S_4\cup S_5=\{a,b,c,d,e,f,h,i,j\}$.
It is easy to see that $|S_4\cup S_2\cup S_6\cup S_5|=10>|S_2\cup S_3\cup S_4\cup S_5|=9$.

For the two examples above, it is easy to check that their $1$-batch and $2$-batch greedy solutions are not unique. For Example~\ref{exampleB1}, the $1$-batch greedy solution $\{S_3, S_4,$
$S_5\}$ is also one solution of the $2$-batch greedy strategy. If we choose $\{S_3, S_4,S_5\}$ instead of $\{S_1, S_5,S_3\}$ as the solution of the $2$-batch greedy strategy, then the $1$-batch greedy strategy has the same performance as the $2$-batch greedy strategy in this case. For Example~\ref{exampleB2}, the $1$-batch greedy solution $\{S_4, S_2,S_6, S_5\}$ is also one solution of the $2$-batch greedy strategy. So we can say that for the $k$-batch greedy strategy, its solution is not unique. However, our harmonic bound under general matroid constraints and exponential bound under uniform matroid constraints are both universal, which means that the harmonic bound holds for any $k$-batch greedy solution under general matroid constraints and the exponential bound holds for any $k$-batch greedy solution under uniform matroid constraints.


\begin{thebibliography}{} 
\bibitem{streeter2008online} Streeter, M.,  Golovin, D.: An online algorithm for maximizing submodular functions. In: Proc. of Adv. Neural Inf. Process. Syst., pp.~1577--1584  (2008)
\bibitem{FeigeVondrak} Feige, U., Vondr{\'a}k, J.: Approximation algorithms for allocation problems: Improving the factor of $1-1/e$. In: Proc. of  47th IEEE Symp. on Foundations of Computer Science, pp.~667--676 (2006)
\bibitem{Fleischer2006}
Fleischer, L., Goemans, M.X., Mirrokni, V.S., Sviridenko, M.: Tight approximation algorithms for maximum general assignment problems. In: Proc. of 17th Annual ACM-SIAM Symp. Discrete Algorithm, pp.~611--620 (2006)
\bibitem{Cohen2006}
Cohen, R., Katzir, L., Raz, D.: An efficient approximation for the generalized assignment problem. Inf. Process. Lett. \textbf{100}(4), 162-166 (2006)
\bibitem{Shmoys1993}
Shmoys, D.B., Tardos, {\'E}.: An approximation algorithm for the generalized assignment problem. Math. Program. \textbf{62}, 461--474  (1993)
\bibitem{Vondrak2011}
Calinescu, G., Chekuri, C., P{\'a}l, M., Vondr{\'a}k, J.: Maximizing a monotone submodular function subject to a matroid constraint. SIAM J. Comput. \textbf{40}(6), 1740--1746 (2011)
\bibitem{Korula2015}
Korula, N., Mirrokni, V., Zadimoghaddam, M.: Online submodular welfare maximization: Greedy beats 1/2 in random order. In: Proc. of 47th Annual Symp. on Theory of Computing,  pp.~889--898 (2015)
\bibitem{Vondrak2008}
Vondr{\'a}k, J.: Optimal approximation for the submodular welfare problem in the value oracle model. In: Proc. of 40th Annual ACM Symp. on Theory of Computing, pp.~67--74 (2008)
\bibitem{K-cover1998}
Hochbaum, D.S., Pathria, A.: Analysis of the greedy approach in problems of maximum $k$-coverage. Nav. Res. Log. \textbf{45}(6), 615--627 (1998)
\bibitem{Feige1998}
Feige, U.: A threshold of $\ln n$ for approximating set cover. J. ACM \textbf{45}(4), 634--652 (1998)
\bibitem{Khuller1999}
Khuller, S., Moss, A., Naor, J.: The budgeted maximum coverage problem. Inf. Process. Lett. \textbf{70}(1), 39--45 (1999)
\bibitem{Fisher1977}
Cornu{\'e}jols, G., Fisher, M.L., Nemhauser, G.L.: Location of bank accounts to optimize float: An analytic study of exact and approximate algorithms. Manag. Sci. \textbf{23}(8), 789--810 (1977)
\bibitem{Location}
Megiddo, N., Zemel, E., Hakimi, S.L.: The maximum coverage location problem. SIAM J. Alg. Disc. Meth. \textbf{4}(2), 253--261 (1983)
\bibitem{Church1974}
Church, R., Velle, C.R.: The maximal covering location problem. Pap. Reg. Sci. \textbf{32}(1), 101--118 (1974)
\bibitem{Pirkul1991}
Pirkul, H., Schilling, D.A.: The maximal covering location problem with capacities on total workload. Manag. Sci. \textbf{37}(2), 233--248 (1991)
\bibitem{LiC12} 
Liu, E., Chong, E.K.P., Scharf, L.L.: Greedy adaptive linear compression in signal-plus-noise models. IEEE Trans.~Inf.~Theory \textbf{60}(4), 2269--2280 (2014)
\bibitem{SensorPlacement}
Krause, A., Singh, A., Guestrin, C.: Near-Optimal sensor placements in Gaussian processes: Theory, efficient algorithms and empirical studies. J. Mach. Learn. Res. \textbf{9}, 235--284 (2008)
\bibitem{Chen2005}
Chen, Y., Chuah, C., Zhao, Q.: Sensor placement for maximizing lifetime per unit cost in wireless sensor networks. In: Proc. of IEEE Military Comm. Conf., pp.~1097--1102 (2005)
\bibitem{Ragi2015}
Ragi, S., Mittelmann, H.D., Chong, E.K.P.: Directional sensor control: Heuristic approaches. IEEE Sens. J. \textbf{15}(1), 374--381 (2015)
\bibitem{nemhauser19781}
Nemhauser, G.L., Wolsey, L.A., Fisher, M.L.: An analysis of approximations for maximizing submodular set functions--{I}.  Math. Program. \textbf{14}(1), 265--294 (1978)
\bibitem{nemhauser1978}
Fisher, M.L., Nemhauser, G.L., Wolsey, L.A.: An analysis of approximations for maximizing submodular set functions--{II}. Math. Program. Stud. \textbf{8}, 73--87 (1978)
\bibitem{hausmann1980}
Hausmann, D., Korte., B., Jenkyns, T.A.: Worst case analysis of greedy type algorithms for independence systems. Math. Program. Stud. \textbf{12}, 120--131 (1980)
\bibitem{conforti1984submodular}
Conforti, M., Cornu{\'e}jols, G.: Submodular set functions, matroids and the greedy algorithm: {T}ight worst-case bounds and some generalizations of the {R}ado-{E}dmonds theorem.  Discrete Appl. Math. \textbf{7}(3), 251--274 (1984)
\bibitem{vondrak2010submodularity}
Vondr{\'a}k, J.: Submodularity and curvature: The optimal algorithm. RIMS Kokyuroku Bessatsu \textbf{B23}, 253--266 (2010)
\bibitem{sviridenko2015}
Sviridenko, M., Vondr{\'a}k, J., Ward, J.: Optimal approximation for submodular and supermodular optimization with bounded curvature. In: Proc. of 26th Annual ACM-SIAM Symp. on Discrete Algorithms, pp. 1134--1148 (2015)
\bibitem{Liu2016}
Liu, Y., Zhang, Z., Chong, E.K.P., Pezeshki, A.: Performance bounds for the $k$-batch greedy strategy in optimization problems with curvature. In: Proc. of 2016 American Control Conf., pp. 7177--7182 (2016)
\bibitem{Chandu2015}
Chandu., D.P.: Big step greedy heuristic for maximum coverage problem. Int. J. Comput. Appl., \textbf{125}(7), 19--24 (2015)
\bibitem{Edmonds}
Edmonds, J.: Submodular functions, matroids, and certain polyhedra. Combinatorial Optimization, \textbf{2570}, 11--26 (2003)
\bibitem{Boros2003}
Boros, E., Elbassioni, K., Gurvich, V., Khachiyan, L.: An inequality for polymatroid functions and its applications. Discrete Appl. Math. \textbf{131}(2), 255--281 (2003)
\bibitem{Tutte}
Tutte, W.T.: Lecture on matroids. J. Res. Natl. Bur. Stand.-Sec. B: Math.  Math.  Phys. \textbf{69B}, 1--47 (1965)
\bibitem{vondrak2010}
Lee, J., Sviridenko, M., Vondr{\'a}k, J.: Submodular maximization over multiple matroids via generalized exchange properties. Math. Oper. Res. \textbf{35}(4), 795--806 (2010)
\bibitem{ZhC13J}
Zhang, Z., Chong, E.K.P., Pezeshki, A., Moran, W.: String submodular functions with curvature constraints. IEEE Trans.~Autom.~Control  \textbf{61}(3), 601--616 (2016)
\bibitem{YJ2015}
Liu, Y.,  Chong, E.K.P.,  Pezeshki, A.: Bounding the greedy strategy in finite-horizon string optimization. In: Proc. of 54th IEEE Conf. Decis. Control, pp. 3900--3905 (2015)
\end{thebibliography}
\end{document}